\documentclass[10pt]{amsart}
\usepackage{amsfonts,amssymb,amscd,amsmath,enumerate,verbatim,calc,amsthm}
\input xy
\xyoption{all}

\headsep=1cm \numberwithin{equation}{section}
\newtheorem{thm}{Theorem}[section]
\newtheorem{cor}[thm]{Corollary}
\newtheorem{lem}[thm]{Lemma}
\newtheorem{prop}[thm]{Proposition}
\newtheorem{defn}[thm]{Definition}

\newtheorem{rem}[thm]{Remark}

\newtheorem{defr}[thm]{Definition and Facts }

\newcommand{\coker}{\mbox{Coker}\,}
\newcommand{\Hom}{\mbox{Hom}\,}
\newcommand{\RHom}{\mbox{RHom}\,}
\newcommand{\Ext}{\mbox{Ext}\,}
\newcommand{\Tor}{\mbox{Tor}\,}
\newcommand{\Spec}{\mbox{Spec}\,}

\newcommand{\Supp}{\mbox{Supp}\,}

\newcommand{\gr}{\mbox{grade}\,}
\newcommand{\Gid}{\mbox{Gid}\,}

\newcommand{\depth}{\mbox{depth}\,}
\renewcommand{\dim}{\mbox{dim}\,}
\renewcommand{\Im}{\mbox{Im}\,}
\newcommand{\cd}{\mbox{cd}\,}
\newcommand{\p}{\mbox{P}\,}

\newcommand{\pd}{\mbox{proj.dim}\,}
\newcommand{\id}{\mbox{id}\,}

\newcommand{\cid}{\mbox{C-id}\,}
\newcommand{\cgid}{\mbox{$G_{C}$-id}\,}
\newcommand{\xid}{\mbox{$\mathcal{X}$-id}\,}

\newcommand{\h}{\mbox{ht}\,}
\newcommand{\im}{\mbox{Im}\,}

\newcommand{\E}{\mbox{E}}

\renewcommand{\H}{\mbox{H}}

\newcommand{\fa}{\mathfrak{a}}

\newcommand{\fm}{\mathfrak{m}}
\newcommand{\fp}{\mathfrak{p}}

\begin{document}
\title[local cohomology modules and Gorenstein injectivity]
 { local cohomology modules and Gorenstein injectivity
with respect to a semidualizing module}

 \author[M.R. Zargar]{Majid Rahro Zargar }
\address{Faculty of mathematical sciences and computer, Tarbiat Moallem
University, 599 Taleghani Avenue, Tehran 15618, Iran, and
School of Mathematics, Institute for Research in Fundamental Sciences (IPM), P.O. Box 19395--5746, Tehran, Iran.}

\email{zargar\_\,m60@tmu.ac.ir}
\email{zargar9077@gmail.com}




\begin{abstract}
Let $(R,\fm)$ be a local ring and let $C$ be a semidualizing $R$-module. In this paper, we are concerned with the $C$-injective and $G_{C}$-injective dimensions of certain local cohomology modules of $R$. Firstly, the injective dimension of $C$ and the above quantities are compared. Secondly, as an application of the above comparisons, a characterization of a dualizing module of $R$ is given. Finally, it is shown that if $R$ is Cohen-Macaulay of dimension $d$ such that $\H_{\fm}^{d}(C)$ is $C$-injective, then $R$ is Gorenstein. This is an answer to the question which was recently raised.
\end{abstract}

\maketitle

\section{introduction}
Throughout this paper, $R$ is a commutative Noetherian ring. It is well known, see for example [6, Corollary 9.5.13], that a local
ring $(R,\fm)$ is Gorenstein if and only if $R$ is a Cohen-Macaulay ring and the
top local cohomology module of $R$, $\H^d_{\fm}(R)$, is isomorphic to $\E_{R}(R/\fm$). As a generalization of this result, in
[12, Theorem 2.5], the author and H. Zakeri showed that $R$ is a Gorenstein local ring if and only if $\id_{R}\H^{\h\fa}_{\fa}(R)<\infty$ for some ideal $\fa$ of $R$ such that $R$ is relative Cohen-Macaulay with respect to $\fa$. In
[17] R. Sazeedeh showed that over a Gorenstein local ring of Krull dimension at
most two the top local cohomology module $\H^{\dim R}_{\fa}(R)$ is a Gorenstein injective
$R$-module for any ideal $\fa$ of $R$. Notice that an injective module is a Gorenstein injective module. In [20]
T. Yoshizawa, as a generalization of Sazeedeh's results, showed that over a complete Cohen-Macaulay local
ring $(R,\fm)$ of Krull dimension $d$ the following conditions are equivalent.
\begin{itemize}
\item[(i)]{$R$ is a Gorenstein ring.}
\item[(ii)]{${\H}^d_{\fm}(R)$ is an injective $R$-module.}
\item[(iii)]{${\H}^d_{\fm}(R)$ is a Gorenstein injective $R$-module.}
\end{itemize}

As a generalization of the above result, it was proved in [14], as the main result, that if $C$ is a semidualizing module over a complete local ring of dimension $d$, then following statements are equivalent.
\begin{itemize}
\item[(i)]{$C$ is a dualizing $R$-module.}
\item[(ii)]{${\H}^d_{\fm}(R)$ is a $C$-injective $R$-module.}
\item[(iii)]{${\H}^d_{\fm}(R)$ is a $G_{C}$-injective $R$-module.}
\item[(iv)]{${\cgid}_{R}{\H}^d_{\fm}(R)<\infty$}.
\end{itemize}
The above result is not true without the Cohen-Macaulay assumption on $R$ (see 3.5). Thus [14, Theorem 3.1] needs correction, nevertheless its proof is clearly valid in the Cohen-Macaulay case.

In this paper, we first prove, in 3.2, that if $R$ is relative Cohen-Macaulay with respect to an ideal $\fa$ of $R$ and $C$ is a semidualizing $R$-module, then $\cid_{R}{\H}^n_{\fa}(R)=\id_{R}C-n=\id_{R}\H_{\fa}^n(C)$, where $n=\h\fa$. Next, in 3.3, we characterize a dualizing module in terms of $G_{C}$-injective dimensions of certain local cohomology modules. Now, as a corollary of 3.2 and 3.3, we establish, in 3.4, the corrected version of [14, Theorem 3.1] without the completeness assumption on $R$. Also, 3.7 provides a generalization of [20, Corollary 2.10] and [12, Corollary 3.14]. Finally, when $R$ is Cohen-Macaulay, we obtain, in 3.8, an answer to the question ``what happens if the top local cohomology module of a semidualizing module $C$ is $C$-injective?" which is stated in [14].

\section{preliminaries}

In this section we recall some definitions and facts which are needed throughout this paper.
\begin{defn}\emph{Following [16, Definition 2.1], let $\mathcal{X}$ be a class of $R$-modules and let $M$ be an $R$-module. An $\mathcal{X}$-\textit{coresolution}
of $M$ is a complex of $R$-modules in $\mathcal{X}$ as follows $$X=0\longrightarrow X_{0}\stackrel{\partial_{0}^X} \longrightarrow X_{-1}\stackrel{\partial_{-1}^X}\longrightarrow \cdots\stackrel{\partial_{n+1}^X}\longrightarrow X_{n}\stackrel{\partial_{n}^X}\longrightarrow X_{n-1}\stackrel{\partial_{n-1}^X}\longrightarrow \cdots$$
such that $\H_{0}(X)\cong M$ and {$\H_{n}(X)=0$} for all $n\leq-1$. The $\mathcal{X}$-\textit{injective dimension} of
$M$ is the quantity
$${\xid_{R}(M)}=\inf\{ \sup \{-n\geq0 | X_{n}\neq0\} ~|~ X~\text{is an $\mathcal{X}$-coresolution of  $M$~}\}.$$
The modules of $\mathcal{X}$-injective dimension zero are precisely the nonzero modules of $\mathcal{X}$ and also $\xid_{R}(0)=-\infty$. }\\

\emph{The following notion of semidualizing modules goes back at least to
Vasconcelos [19], but was rediscovered by others. The reader is referred to [15] for more details about semidualizing modules.}

\end{defn}
\begin{defn}\emph{A finitely generated $R$-module $C$ is called \textit{semidualizing} if
the natural homomorphism $R\rightarrow \Hom_{R}(C,C)$ is an isomorphism and
$\Ext_{R}^{i}(C,C)=0$ for all $i\geq1$. An $R$-module $D$ is said to be a \textit{dualizing} $R$-module if it is semidualizing and
has finite injective dimension. For a semidualizing $R$-module $C$, the class of $C$-\textit{injective} modules is defined as
$$\mathcal{I}_{C}(R)=\{~\Hom_{R}(C,I) | ~~ I  ~~ \text{is an injective $R$-module}\}.$$
It will be convenient for us to denote $\mathcal{I}_{C}(R)$-$\id_{R} M$, which is defined in 2.1, by $\cid_{R}M$.
Notice that when $C=R$ these notions recover the concepts of injective module and injective dimension. }
\end{defn}

Based on the work of E.E. Enochs and O.M.G. Jenda [6], the following notions were introduced
and studied by H. Holm and P. J${\o}$rgensen [9].
\begin{defn}\emph{Let $C$ be a semidualizing $R$-module. A complete $\mathcal{I}_{C}\mathcal{I}$-\textit{resolution} is a complex $Y$ of $R$-modules such that:
\begin{itemize}
\item[(i)]{$Y$ is exact and {$\Hom_{R}(I,Y)$} is exact for each $I\in \mathcal{I}_{C}(R)$.}
\item[(ii)]{ $Y_{i}\in \mathcal{I}_{C}(R)$ for all $i>0$ and $Y_{i}$ is injective for all $i\leq0$.}
\end{itemize}
An $R$-module $M$ is called $G_{C}$-\textit{injective} if there exists a complete $\mathcal{I}_{C}\mathcal{I}$-resolution $Y$
such that $M\cong \ker(\partial_{0}^{Y})$. In this case $Y$ is a complete $\mathcal{I}_{C}\mathcal{I}$-resolution of $M$. The class of $G_{C}$-injective $R$-modules is denoted by $\mathcal{GI_{C}}(R)$. For convenience the $\mathcal{GI_{C}}(R)$-\textit{injective dimension}, $\mathcal{GI_{C}}(R)$-$\id_{R}M$, of $M$ which is defined as in 2.1 is denoted by $\cgid_{R}(M)$.}

\emph{Note that when $C=R$ these notions recover the concepts of Gorenstein injective module and Gorenstein injective dimension which were introduced in [6].
}
\end{defn}

\begin{defn} \emph{We say that a finitely generated $R$-module $M$ is \textit{relative Cohen–
Macaulay with respect to $\fa$} if there is precisely one non-vanishing local cohomology module of $M$ with
respect to $\fa$. Clearly this is the case if and only if $\gr(\fa,M)=\cd(\fa,M)$, where $\cd(\fa,M)$ is the largest integer $i$ for which $\H_{\fa}^i(M)\neq0$. Observe that the notion of relative Cohen-Macaulay module is connected with the notion of cohomologicaly complete intersection ideal which has been studied in [7].
}
\end{defn}
\begin{rem}\emph{Let $M$ be a relative Cohen-Macaulay module with respect to $\fa$ and let $\cd(\fa,M)=n$. Then, in view of [2, Theorems 6.1.4, 4.2.1, 4.3.2], it is easy to see that $\Supp\H^{n}_{\fa}(M)=\Supp({M}/{\fa M})$ and $\h_{M}\fa=\gr(\fa,M)$, where
$\h_{M}\fa =\inf\{\ \dim_{R_{\fp}}M_{\fp} |~ \fp\in\Supp(M/\fa M) ~\}$.}
\end{rem}
Next, we recall some elementary results about the trivial extension of a ring by a module.
\begin{defr}\emph{Let $C$ be an $R$-module. Then the direct sum $R\oplus C$ has the structure of a commutative ring with respect to multiplication defined by
$$(a,c)(a',c')=(aa', ac'+a'c),$$}
\emph{for all $(a,c),(a', c')\in R\oplus C$. This ring is called \textit{trivial extension} of $R$ by $C$ and is denoted by $R\ltimes C$. The following properties of $R\ltimes C$ are needed in section 3. }
\begin{itemize}
\item[(i)]{There are natural ring homomorphisms $R\rightleftarrows R\ltimes C$ which enable us to consider $R$-modules as $R\ltimes C$-modules, and vice versa.}
\item[(ii)]\emph{{For any ideal $\fa$ of $R$, $\fa\oplus C$ is an ideal of $R\ltimes C$.}
\item[(iii)]{$(R\ltimes C, \fm\oplus C)$ is a Noetherian local ring whenever $(R,\fm)$ is a Noetherian local ring and $C$ is a finitely generated $R$-module. Also, in this case, $\dim R=\dim R\ltimes C$}.}
\item[(v)]\emph{For any $R$-module $M$ we have  $\cgid_{R}M =\Gid_{R\ltimes C}M$ (see [9, Theorem 2.16])}.

\end{itemize}
\end{defr}
The classes defined next are collectively known as Foxby classes. The reader is referred to [1], [15] and [18] for some basic results about those classes.
\begin{defn}\emph{Let $C$ be a semidualizing $R$-module. The \textit{Bass class} with respect to $C$ is the class $\mathcal{B_{C}}(R)$ of $R$--modules such that: }
\begin{itemize}
\item[(i)]{\emph{$\Ext_{R}^i(C,M)=0=\Tor^{R}_{i}(C,\Hom_{R}(C,M))$ for all $i\geq1$}, and}
\item[(ii)]{\emph{the natural evaluation map $C\otimes_{R}\Hom_{R}(C,M)\rightarrow M$ is an isomorphism}.}
\end{itemize}
\emph{Dually, the \textit{Auslander class} with respect to $C$, denoted $\mathcal{A}_{C}(R)$, consists of
all $R$-modules $M$ such that:}
\begin{itemize}
\item[(i)]{\emph{$\Tor^{R}_{i}(C,M)=0=\Ext_{R}^{i}(C,C\otimes_{R}M)$ for all $i\geq1$}, and}
\item[(ii)]{\emph{the natural map $M\rightarrow \Hom_{R}(C,C\otimes_{R}M)$ is an isomorphism}.}
\end{itemize}
\end{defn}
\section{main results}
The starting point of this section is the next lemma, which is assistance in the proof
of Theorem 3.2.
\begin{lem}Let $C$ be a semidualizing $R$-module and let
$$0\rightarrow M'\rightarrow M\rightarrow M''\rightarrow 0$$
be an exact sequence of $R$-modules and $R$-homomorphisms such that $M'$ and $M$ are $C$-injective. Then $M''$ is $C$-injective.
\begin{proof}First, in view of [18, Proposition 3.4] we see that $\cid_{R}M''$ is finite. Therefore, by [18, Corollary 2.9], $M''\in{\mathcal{A}}_{C}(R)$ and hence $\Tor^{R}_{1}(C,M'')=0$. Thus the sequence $$0\rightarrow C\otimes_{R}M'\rightarrow C\otimes_{R}M\rightarrow C\otimes_{R}M''\rightarrow 0$$
is exact. Now, one can use [18, Theorem 2.15] to complete the proof.
\end{proof}
\end{lem}

\begin{thm}Suppose that $(R,\fm)$ is a local ring and that $C$ is a semidualizing $R$-module. Suppose that $M$ is relative Cohen-Macaulay with respect to $\fa$ and that \emph{$\h_{M}\fa=n.$} Then the following statements hold.
\begin{itemize}
\item[(i)]{\emph{$\cid_{R}{\H}^n_{\fa}(M)\leq\cid_{R}M-n$}}
\item[(ii)]{In the case where $M=R$, \emph{$C$} is a dualizing $R$--module if and only if \emph{$\cid_{R}{\H}^n_{\fa}(R)$ is finite.}}
Furthermore, we have the equality \emph{$\cid_{R}{\H}^n_{\fa}(R)=\id_{R}C-n=\id_{R}\H_{\fa}^{n}(C)$}.
\end{itemize}
\begin{proof}First we present a $C$-injective resolution approach for calculation of local cohomology modules. To this end, let $\fa$ be an ideal of $R$, $N$ be a finitely generated $R$-module and let $I$ be an injective $R$-module. Then we have the following natural isomorphisms
{ \[\begin{array}{rl}
\Gamma_{\fa}(\Hom_{R}(N,I))&\cong\underset{n\in\mathbb{N}}\varinjlim \Hom_{R}(R/\fa^n,\Hom_{R}(N,I))\\
&\cong\underset{n\in\mathbb{N}}\varinjlim \Hom_{R}(R/\fa^n\otimes_{R}N,I)\\
&\cong\underset{n\in\mathbb{N}}\varinjlim \Hom_{R}(N,\Hom_{R}(R/\fa^n,I))\\
&\cong\Hom_{R}(N,\underset{n\in\mathbb{N}}\varinjlim \Hom_{R}(R/\fa^n,I))\\
&\cong\Hom_{R}(N,\Gamma_{\fa}(I)).
\end{array}\]}
Therefore,  $\Gamma_{\fa}(N')\in\mathcal{I}_{C}(R)$ for all $N'\in\mathcal{I}_{C}(R)$, because $\Gamma_{\fa}(I)$ is injective. Also, by using a finite free resolution for $N$ and the above isomorphism, we can deduce that $\H_{\fa}^{i}(\Hom_{R}(N,I))=0$ for all $i\geq1$. Hence, $C$-injective modules are $\Gamma_{\fa}$-acyclic. Therefore in view of [2, Exercise 4.1.2], for an $R$-module $N$, its local cohomology modules with respect to $\fa$ can be calculated by means of a $C$-injective resolution for $N$. Now, we prove the assertion.

(i) Let $\cid_{R}(M)=d$ and let
\begin{equation}
0\rightarrow M \rightarrow \Hom_{R}(C,\E^{0})\rightarrow \cdots\stackrel{\alpha}\rightarrow \Hom_{R}(C,\E^{n})\stackrel{\beta}\rightarrow
\cdots\rightarrow \Hom_{R}(C,\E^{d})\rightarrow0
\end{equation}
be a $C$-injective resolution of $M$. Then, by applying the functor $\Gamma_{\fa}(-)$ on (3.1), we may use the above discussions in conjunction with our assumption on local cohomology module of $M$ to obtain the exact sequences
\begin{equation}0\rightarrow \Gamma_{\fa}(M) \rightarrow\Hom_{R}(C,\Gamma_{\fa}(\E^{0}))\rightarrow \cdots\rightarrow\Hom_{R}(C,\Gamma_{\fa}(\E^{n}))\rightarrow\coker{\Gamma_{\fa}(\alpha)}\rightarrow0
\end{equation}
and
\begin{equation}0\rightarrow \Im{\Gamma_{\fa}({\beta})}\hookrightarrow\Hom_{R}(C,\Gamma_{\fa}(\E^{n+1}))\rightarrow \cdots\rightarrow\Hom_{R}(C,\Gamma_{\fa}(\E^{d}))\rightarrow0.
\end{equation}
Now, if $n=0$ the result is clear. So suppose that $n>0$. Then, by assumption, $\Gamma_{\fa}(M)=0$. Therefore, by using the exact sequence (3.2) and
Lemma 3.1, we see that $\coker{\Gamma_{\fa}(\alpha)}$ is $C$-injective. Notice that $\H_{\fa}^n(M)=\frac{\ker\Gamma_{\fa}(\beta)}{\im\Gamma_{\fa}(\alpha)}$.
Therefore, patching the long exact sequence (3.3) together with the exact sequence$$0\longrightarrow\H_{\fa}^n(M)\longrightarrow\coker{\Gamma_{\fa}(\alpha)}\longrightarrow\Im{\Gamma_{\fa}({\beta})}
\longrightarrow0,$$
gives the following long exact sequence
$$0\rightarrow\H_{\fa}^n(M)\rightarrow\coker{\Gamma_{\fa}(\alpha)}\rightarrow\Hom_{R}(C,\Gamma_{\fa}(\E^{n+1}))\rightarrow \cdots\rightarrow\Hom_{R}(C,\Gamma_{\fa}(\E^{d}))\rightarrow0$$
Hence, $\cid_{R}{\H}^n_{\fa}(M)\leq\cid_{R}M-n$.

(ii) Let $R$ be relative Cohen-Macaulay with respect to $\fa$. First we notice that, $\Supp_{R}(C)=\Spec(R)$. Therefore, in view of [5, Theorem 2.2], $\cd(\fa,R)=\cd(\fa,C)$ and, by [15, Theorem 2.2.6(c)], $\gr(\fa,R)=\gr(\fa,C)$. Hence, in view of 2.5, $C$ is relative Cohen-Macaulay with respect to $\fa$. Since, by [2, Theorem 3.4.10], for each $i\geq0$, the local cohomology functor $\H_{\fa}^i(-)$ commutes with direct limits and any $R$-module can be viewed as the direct limit of its finitely generated submodules, one can use [5, Theorem 2.2] to see that the functor $\H_{\fa}^{n}(-)$ is right exact. Therefore, in view of [2, Exercise 6.1.9], we have $\H_{\fa}^n(R)\otimes_{R}C\cong\H_{\fa}^n(C)$. Hence, by [18, Theorem 2.11], $\cid_{R}{\H}^n_{\fa}(R)=\id_{R}{\H}^n_{\fa}(C)$. Now, one can use [12, Theorem 2.5] to complete the proof.
 \end{proof}
\end{thm}
The next proposition is a generalization of [12, Corollary 3.10].
\begin{prop} Let $(R,\fm)$ be a Cohen-Macaulay local ring which has a dualizing
module and let $C$ be a semidualizing $R$-module. Then the following conditions are equivalent.
\begin{itemize}
\item[(i)]{$C$ is a dualizing $R$-module.}
\item[(ii)]{$\emph{\cgid}_{R}\emph{\H}^n_{\fa}(R)<\infty$ for some ideal $\fa$ of $R$ such that $R$ is relative Cohen-Macaulay with respect to $\fa$ and that $\emph{\h}_{R}\fa =n$}.
\end{itemize}
\begin{proof} Suppose that $\fa$ is an ideal of $R$ such that $R$ is relative Cohen-Macaulay with respect to $\fa$. and set $\h_{R}\fa=n$. Then the implication (i)$\Rightarrow$(ii) follows from 3.2(ii).

(ii)$\Rightarrow$(i) Let $\cgid_{R}\H^n_{\fa}(R)<\infty$. Then, by 2.6(v), we have $\Gid_{R\ltimes C}\H^n_{\fa}(R)<\infty$. Now, in view of the Independence Theorem [2, Theorem 4.2.1], we have $\H_{\fa}^i(R)\cong\H_{\fa\oplus C}^i(R)$ for all $i$. Hence, $R$ is relative Cohen-Macaulay with respect to $\fa\oplus C$ as an $R\ltimes C$-module. On the other hand, in view of [3, Exercise 1.2.26] and [15, Theorem 2.2.6], we deduce that $R\ltimes C$ is a Cohen-Macaulay local ring. Also, by the remark before [9, Lemma 4.5], $R\ltimes C$ has dualizing complex $D_{R\ltimes C}=\RHom_{R}(R\ltimes C, D_{R})$, where $D_{R}$ is a dualizing complex of $R$. Since $\Gid_{R\ltimes C}\H^n_{\fa \oplus C}(R)<\infty$ and $R$ is Cohen-Macaulay as an $R\ltimes C$-module, in view of [12, Theorem 3.8], we have $\Gid_{R\ltimes C}R<\infty$. Therefore, by [10, Proposition 4.5], $\id_{R}C\leq\Gid_{R\ltimes C}R<\infty$. Thus $C$ is a dualizing $R$--module.
\end{proof}
\end{prop}
Next, we single out a certain case of Theorem 3.2(ii) and Proposition 3.3 in which $\fa=\fm$. Here, of course, we are not assuming that $R$ has a dualizing module.
\begin{cor} Let $(R,\fm)$ be a Cohen-Macaulay local ring with dimension $n$ and let $C$ be a semidualizing $R$-module. Then the following statements are equivalent.
\begin{itemize}
\item[(i)]{$C$ is a dualizing $R$-module.}
\item[(ii)]{$\emph{\H}^n_{\fm}(R)$ is a $C$-injective $R$-module.}
\item[(iii)]{$\emph{\H}^n_{\fm}(R)$ is a $G_{C}$-injective $R$-module.}
\item[(iv)]{$\emph{\cgid}_{R}\emph{\H}^n_{\fm}(R)<\infty$}.
\end{itemize}
\begin{proof}First notice that $R$ is relative Cohen-Macaulay with respect to $\fm$. Thus, the implications (i)$\Leftrightarrow$(ii) follows from 3.2 and [3, Theorem 3.1.17]. Also, the implication (ii)$\Rightarrow$(iii) holds by [9, Example 2.8].

(iii)$\Rightarrow$(iv) is clear.

(iv)$\Rightarrow$(i) By the same arguments as in the proof of 3.3(ii)$\Rightarrow$(i) one can see that $R$ and $R\ltimes C$ are Cohen-Macaulay modules over $R\ltimes C$ and that $\Gid_{R\ltimes C}\H^n_{\fm \oplus C}(R)<\infty$. Hence by [12, Corollary 3.9] and [10, Proposition 4.5] we have $\id_{R}C\leq\Gid_{R\ltimes C}R<\infty$.
\end{proof}
\end{cor}

\begin{rem}\emph{As a main result, it has been proved in [14, Theorem 3.1] that if $(R,\fm)$ is a complete local ring with $\dim R=n$ and $C$ is a semidualizing $R$-module, then the statements (i),(ii),(iii) and (iv) of 3.4 are equivalent. This result is not true without the Cohen-Macaulay assumption on $R$. Thus [14, Theorem 3.1] needs correction, nevertheless its proof is clearly valid in the Cohen-Macaulay case. Indeed, M. Hermann and N. V. Trung, in [8], present a Buchsbaum ring $(R,\fm,k)$ with $\dim R =3$ which is not Gorenstein, but $\H_{\fm}^3(R)\cong E_{R}(k)$. However, 3.4 recovers the corrected version of [14, Theorem 3.1].}
\end{rem}

The New Intersection Theorem implies that if a local ring admits a finitely generated module of finite injective dimension, then the ring is Cohen-Macaulay. (This was formerly known as Bass' Conjecture.) For the proof of this result the reader is referred to [11] and [13]. In the next proposition, which is a generalization of [14, Proposition 3.3], we shall use this result.
\begin{prop}Let $(R,\fm)$ be a local ring which has a dualizing module $D$ and let $C$ be a semidualizing $R$-module.
Set \emph{$C':=\Hom_{R}(C,D)$}. Then the following statements are equivalent.
\begin{itemize}
\item[(i)]{\emph{$C'\cong R$}.}
\item[(ii)]{\emph{$\H_{\fa}^{\h_{R}\fa}(R)\in\mathcal{B_{C'}}(R)$} for any ideal $\fa$ of $R$ such that $R$ is relative Cohen-Macaulay with respect to $\fa$.}
\item[(iii)]{\emph{$\H_{\fa}^{\h_{R}\fa}(R)\in\mathcal{B_{C'}}(R)$} for some ideal $\fa$ of $R$ such that $R$ is relative Cohen-Macaulay with respect to $\fa$.}
\end{itemize}
\begin{proof} First we notice that, since $R$ has a dualizing module, there exists a finitely generated $R$-module with finite injective dimension. Therefore, by the New Intersection Theorem, $R$ is Cohen-Macaulay. Now, since $\mathcal{B}_{R}(R)$ is precisely the category of $R$-modules, the implication (i)$\Rightarrow$(ii) holds. Also, since $R$ is relative Cohen-Macaulay with respect to $\fm$, the implication (ii)$\Rightarrow$(iii) is true obviously.

(iii)$\Rightarrow$(i) Suppose that $\fa$ is an ideal of $R$ such that $R$ is relative Cohen-Macaulay with respect to $\fa$ and that $\H_{\fa}^{\h_{R}\fa}(R)\in\mathcal{B_{C'}}(R)$. Then, in view of [9, Theorem 4.6], $\cgid_{R}\H^{\h_{R}{\fa}}_{\fa}(R)<\infty$. Hence, one can use 3.3 to complete the proof.
\end{proof}
\end{prop}

In [20, Corollary 2.10] it was shown that if ($R,\fm)$ is a
complete Gorenstein local ring with $\dim R=d\leq2$ and $M$ is an $R$-module, then the top local
cohomology module is Gorenstein injective for all proper ideals $\fa$ of $R$. Recently, in [12, Corollary 3.14], it was shown that the converse of this result is also true. Our next theorem is a generalization of this result.

\begin{thm}Let $(R,\fm)$ be a local ring with  $d=\dim R\leq2$, $C$ be a semidualizing $R$-module. Then the following statements are equivalent.
\begin{itemize}
\item[(i)]{$C$ is a dualizing $R$-module.}
\item[(ii)]{\emph{$\H^d_{\fm}(R)$} is a ${C}$-injective $R$-module.}
\item[(iii)]{\emph{$\H^d_{\fm}(R)$} is a $G_{C}$-injective $R$-module.}
\item[(iv)]{ \emph{$\H^d_{\fa}(M)$} is $G_{C}$-injective for all finitely generated $R$-modules $M$ and for all ideals $\fa$.}
\end{itemize}
\begin{proof}(i)$\Rightarrow$(ii) Since $C$ has finite injective dimension, by the New Intersection Theorem, $R$ is a Cohen-Macaulay ring. Hence, one can use 3.2(ii) and [3, Theorem 3.1.17] to see that $\H^d_{\fm}(R)$ is $C$-injective. The implication (ii)$\Rightarrow$(iii) is clear by [9, Example 2.8].

(iii)$\Rightarrow$(i) Suppose that $\H^d_{\fm}(R)$ is $G_{C}$-injective. Then, by 2.6(v) and [2, Theorem 4.2.1], $\H^d_{\fm\oplus C}(R)$ is Gorenstein
injective as an $R\ltimes C$-module. Now, since $\dim R\ltimes C\leq 2$, in view of [12, Proposition 3.12] we have $${\Gid}_{R\ltimes C}{\H}^n_{\fm\oplus C}(R)=\depth(R\ltimes C)-\dim_{(R\ltimes C)}R=\depth R-\dim R.$$ So, $R$ is Cohen-Macaulay. Therefore, by 3.4, $C$ is a dualizing $R$--module.

(i)$\Rightarrow$(iv) Since $C$ is a dualizing $R$-module, in view of [10, Lemma 3.5], the local ring $R\ltimes C$ is Gorenstein. Now, let $\fa$ be an ideal of $R$ and let $M$ be a finitely generated $R$-module. We can assume that $\H_{\fa}^{d}(M)\neq0$ because the conclusion is easy if $\H_{\fa}^{d}(M)=0$. Then, in view of 2.6(v), we have ${\cgid}_{R}{\H}^d_{\fa}(M)={\Gid}_{R\ltimes C}{\H}^d_{\fa}(M)$. Therefore, by [2, Theorem 4.2.1] and [4, Theorem 3.14], we see that ${\cgid}_{R}{\H}^d_{\fa}(M)={\Gid}_{R\ltimes C}{\H}^d_{\fa\oplus C}(M)$ is finite. Now, one can use [12, Corollary 3.14] to complete the proof. The implication (iv)$\Rightarrow$(iii) is clear.
\end{proof}
\end{thm}
Let $(R,\fm)$ be a local ring and let $C$ be a semidualizing $R$-module. The following question is stated in [14]. What happens if the top local cohomolgy module of $C$ is $C$-injective? Next we provide an answer to this question.
\begin{thm}Let $(R,\fm)$ be a Cohen-Macaulay local ring with $\dim R=d$ and let $C$ be a semidualizing $R$-module. Then the following conditions are equivalent.
\begin{itemize}
\item[(i)]{$R$ is Gorenstein.}
\item[(ii)]{\emph{$\cid_{R}\H^d_{\fm}(C)<\infty$}.}
\item[(iii)]{\emph{$\H_{\fm}^d(C)$} is $C$-injective.}
\end{itemize}

\begin{proof}(iii)$\Rightarrow$(i). First notice that, by [15, Proposition 2.2.1], $C$ is a semidualizing $R$-module if and only if $C\otimes_{R}\hat R$ is a semidualizing $\hat R$-module. Now, since, by [2, Exercise 6.1.9], the functor $\H_{\fm}^d(-)$ is right exact, we see, in view of [18, Theorem 2.11], that $\H_{\fm}^d(C)$ is $C$-injective if and only if $\H_{\fm}^d(C\otimes_{R}C)$ is injective. Therefore, one can use [2, Theorem 4.3.2] to deduce that $\H_{\fm}^d(C)$ is a $C$-injective $R$-module if and only if $\H_{\fm\hat R}^d(C\otimes_{R}\hat R)$ is a $(C\otimes_{R}\hat R)$-injective $\hat R$-module. Thus, we may assume that $R$ is complete.

If $d=0$, then $\Gamma_{\fm}(C)=C$. Hence, by assumption, $C$ is $C$-injective: and so, by [16, Lemma 2.11], $R$ is Gorenstein. Let $d>0$. Since $R$ is a complete Cohen-Macauly local ring, it has a canonical module. Therefore, in view of [6, Proposition 9.5.22], $\pd_{R}\H_{\fm}^d(R)=d$. Hence, by [16, Fact 1.6 ], $\H_{\fm}^d(R)\in \mathcal{A}_{C}(R)$; so that $\Tor^{R}_{i}(C,\H_{\fm}^d(R))=0$ for all $i\geq0$. Let \begin{equation}0\longrightarrow {\p}_{d} \longrightarrow \cdots \longrightarrow {\p}_{1}\longrightarrow {\p}_{0}\longrightarrow \H_{\fm}^d(R)\longrightarrow0 \end{equation} be a projective resolution for $\H_{\fm}^d(R)$. Then, by applying the functor $-\otimes_{R}C$ on (3.4) we obtain the exact sequence \begin{equation}0\longrightarrow {\p}_{d}\otimes_{R}C \longrightarrow \cdots \longrightarrow {\p}_{1}\otimes_{R}C\longrightarrow {\p}_{0}\otimes_{R}C\longrightarrow \H_{\fm}^d(R)\otimes_{R}C\longrightarrow0, \end{equation} which is a $C$-projective resolution of $\H_{\fm}^d(R)\otimes_{R}C\cong\H_{\fm}^d(C)$. Now, by [12, Proposition 2.1], $\Ext^i_R(R/\fm,\H^{d}_\fm(C))\cong\Ext^{i+d}_R(R/\fm,C)$ for all $i\geq0$. Therefore, one can use [16, Theorem 2.14] to complete the proof.

(ii)$\Rightarrow$(iii). Assume that $\cid_{R}\H_{\fm}^d(C)<\infty$. Then, $\cgid_{R}\H_{\fm}^d(C)<\infty$. Hence, by similar arguments as in the proof of Proposition 3.3(ii)$\Rightarrow$(i) we have
$\cgid_{R}\H_{\fm}^d(C)=\Gid_{R\ltimes C}\H^d_{\fm \oplus C}(C)$. Next, by [12, Proposition 3.12], $${\Gid}_{R\ltimes C}{\H}^d_{\fm\oplus C}(C)=\depth(R\ltimes C)-\dim_{(R\ltimes C)}C=\depth R-\dim C=0.$$
Therefore, $\H_{\fm}^d(C)$ is $G_{C}$-injective. On the other hand, by our assumption and [18, Corollary 2.9](b), $\H_{\fm}^d(C)\in \mathcal{A}_{C}(R)$. Thus, by the proof of [9, Theorem 4.2], we see that $\H_{\fm}^d(C)\otimes_{R}C\cong\H_{\fm}^d(C\otimes_{R}C)$ is a Gorenstein injective $R$-module. Now, one can use [4, Proposition 3.10] and [18, Theorem 2.11] to deduce that $\H_{\fm}^d(C)$ is $C$-injective. The implication (i)$\Rightarrow$(ii) follows from 3.4(i)$\Rightarrow$(ii) since, when $R$ is Gorenstein, the module $C\cong R$ is dualizing.
\end{proof}
\end{thm}

{$\mathbf{Acknowledgements}.$} I am very grateful to Professor Hossein Zakeri for his kind comments and assistance in the preparation of this article. Also, the author is grateful to the referee for careful reading of the manuscript, for suggesting several improvements of the manuscript and for imposing related questions.

\end{document}